\documentclass[11pt]{article}

\usepackage{tikz}

\usepackage{bm,amssymb,amsthm,amsfonts,amsmath,latexsym,cite,color, psfrag,graphicx,ifpdf,pgfkeys,pgfopts,xcolor,extarrows}

\newtheorem{theo}{Theorem}[section]

\newtheorem{lemma} [theo]{Lemma}

\allowdisplaybreaks

\makeatletter \@addtoreset{equation}{section}
\@addtoreset{theo}{}\makeatother

\usepackage{hyperref}

\def\qed{\hfill \rule{4pt}{7pt}}
\def\pf{\noindent {\it Proof.} }

\def\G{{\mathfrak{G}}}\def\S{{\mathfrak{S}}}

\begin{document}

\begin{center}
{\large\bf A Note on Specializations of Grothendieck Polynomials}
\end{center}

\begin{center}

{\small Neil J.Y. Fan$^1$ and Peter L. Guo$^2$}

\vskip 4mm

$^1$Department of Mathematics\\
Sichuan University, Chengdu, Sichuan 610064, P.R. China

$^{2}$Center for Combinatorics, LPMC\\
Nankai University,
Tianjin 300071,
P.R. China \\[3mm]

\vskip 4mm

$^1$fan@scu.edu.cn,
$^2$lguo@nankai.edu.cn
\end{center}

\begin{abstract}
Buch and Rim\'{a}nyi proved  a formula for a specialization of double  Grothendieck polynomials  based on the Yang-Baxter equation related to the degenerate Hecke algebra.
A geometric proof was found by Yong and Woo    by constructing a Gr\"{o}bner basis for the Kazhdan-Lusztig ideals. In this note, we give an elementary proof for
 this formula by using only   divided difference operators.
\end{abstract}

\section{Introduction}

Let $S_n$ denote the symmetric group of permutations
of $\{1,2,\ldots,n\}$. For a permutation
$w\in S_n$, the double Grothendieck polynomial $\G_w(x;y)$
 introduced by Lascoux and Sch\"utzenberger \cite{LS}  is
the polynomial representative of the class of the Schubert variety for $w$ in the equivariant $K$-theory  of the flag manifold. Write a permutation $v\in S_n$ in  one-line notation, that is,
 write $v=v(1)v(2)\cdots v(n)$. The specialization
\begin{equation}\label{BR-formula}
\G_w(y_v;y):=\G_w(y_{v(1)},\ldots,y_{v(n)};y)
\end{equation}
of $\G_w(x;y)$ obtained by replacing $x_i$ with $y_{v(i)}$
gives the restriction of this class to the fixed point corresponding to $v$. Buch and Rim\'{a}nyi \cite{BR} proved a   formula for  $\G_w(y_v;y)$ based on the Yang-Baxter equation related to the degenerate Hecke algebra. Buch and Rim\'{a}nyi \cite{BR}
 also pointed out various important applications  of this formula. By constructing a Gr\"{o}bner basis for the Kazhdan-Lusztig ideals, Yong and Woo \cite{YW} found   a geometric explanation for
the Buch-Rim\'{a}nyi formula.

In this note, we give an elementary proof  of the Buch-Rim\'{a}nyi
formula  by using only   divided difference operators.
As observed by Buch and Rim\'{a}nyi \cite[Corollary 2.3]{BR}, the
classical pipe dream (or, RC-graph) formula of
$\G_w(x;y)$ (see for example \cite[Corollary 5.4]{KM}, \cite[Theorem 6.3]{LRS}) can be directly obtained from the specialization
$\G_w(y_v;y)$. Hence our approach implies that the
pipe dream formula for double Grothendieck polynomials can be
derived  directly from   divided difference operators.

\section{The Buch-Rim\'{a}nyi formula}

Fix a nonnegative integer $n$. For $1\leq i<j\leq n$, let $t_{ij}$ denote the transposition $(i,j)$
in $S_n$.  So, if $w\in S_n$, then $wt_{ij}$ is the permutation obtained from $w$ by interchanging $w(i)$ and
$w(j)$, while $t_{ij} w$ is obtained from $w$ by interchanging the values $i$ and $j$.
For example, for $w=2143$, we have $wt_{13}=4123$ and $t_{13}w=2341$.
Write $s_i$ for the adjacent transposition $(i, i+1)$.
Each
permutation can be written as a product of adjacent transpositions. The length $\ell(w)$ of
a permutation $w$ is the minimum $k$ such that $w=s_{i_1}s_{i_2}\cdots s_{i_k}$, and in this case,
$(s_{i_1},s_{i_2},\ldots, s_{i_k})$ is called a reduced word of $w$.
It is well known  that the length $\ell(w)$  is equal to the number of
pairs $(i,j)$ such that $i<j$ and $w(i)>w(j)$:
\[\ell(w)=\#\{(i,j)\colon 1\leq i<j\leq n,\  w(i)>w(j)\}.\]
Hence, it is clear that $\ell(ws_i)=\ell(w)+1$ if and only if $w(i)<w(i+1)$, while
 $\ell(ws_i)=\ell(w)-1$ if and only if $w(i)>w(i+1)$.

Let $\mathbb{Z}[x^{\pm},y^{\pm}]$ denote the ring of Laurent polynomials in the $2n$ commuting indeterminates $x_1,\ldots,x_n,y_1,\ldots,y_n$.  For a Laurent polynomial $f(x,y)\in \mathbb{Z}[x^{\pm},y^{\pm}]$,
the  divided difference operator $\partial_i$ acting on $f(x,y)$ is defined by
\[\partial_i f=(f-s_if)/(x_i-x_{i+1}),\]
where  $s_if$ is obtained from $f$ by interchanging $x_i$ and $x_{i+1}$.
It is easy to check that $\partial_if$ is still a Laurent polynomial.
Let $w_0=n\cdots 21$ be the longest permutation in $S_n$. Set
\begin{equation}\label{DDD}
\G_{w_0}(x;y)=\prod_{i+j\le n}\left(1-\frac{y_j}{x_i}\right).
\end{equation}
For $w\neq w_0$, choose an adjacent transposition  $s_i$ such that $\ell(ws_i)=\ell(w)+1$.
Let $\pi_i=\partial_i x_i$ and define
\begin{align}\label{def}
\G_{w}(x;y)&=\pi_i \G_{ws_i}(x;y)=\frac{x_i\G_{ws_i}(x;y)-x_{i+1}\G_{ws_i}(\ldots,x_{i+1},x_i,\ldots; y)}{x_i-x_{i+1}}.
\end{align}
The above definition   is independent of the choice of
 $s_i$ since the operators $\pi_i$ satisfy the
Coxeter  relations: $\pi_i \pi_j=\pi_j \pi_i$ for $|i-j|>1$, and
$\pi_i \pi_{i+1} \pi_i= \pi_{i+1}\pi_i  \pi_{i+1}$, see for example \cite[(2.14)]{Ma}.

We remark  that there are other equivalent definitions for
double Grothendieck polynomials. The definition adopted here implies that
$\G_w(x;y)$ are    Laurent polynomials. The double Grothendieck  polynomials  $\mathfrak{L}_w^{(-1)}(y;x)$ defined  in  \cite{FK1} are legitimate  polynomials, which  can be obtained from $\G_w(x;y)$ by replacing $x_i$ and $y_i$ respectively  with $\frac{1}{1-x_i}$ and $1-y_i$. It should also be
 noticed that $\G_w(x^{-1}; y^{-1})$ are the double Grothendieck polynomials used   in \cite{KM1},
  and $\G_w(x^{-1};y)$ are the double Grothendieck polynomials appearing   in \cite{KM}.
It is worth mentioning that the
double Schubert polynomial  $\S_w(x;y)$ is  the lowest degree homogeneous component
 of $\mathfrak{L}_w^{(-1)}(y;x)$,  see \cite{BB,BJS,FK2,FS,LLS} for  combinatorial constructions of
 Schubert polynomials.

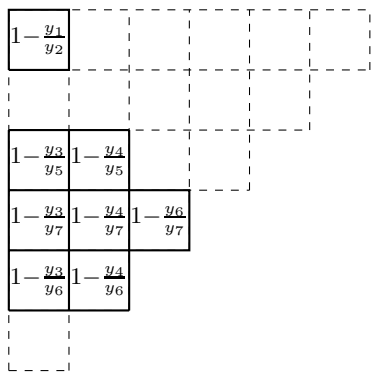
\begin{figure}[h]
\begin{center}
\begin{tikzpicture}(100,140)

\draw [dashed] (0mm,0mm) -- (8mm,0mm);
\draw [dashed] (0mm,8mm) -- (16mm,8mm);
\draw [dashed] (0mm,16mm) -- (24mm,16mm);
\draw [dashed] (0mm,24mm) -- (32mm,24mm);
\draw [dashed] (0mm,32mm) -- (40mm,32mm);
\draw [dashed] (0mm,40mm) -- (48mm,40mm);
\draw [dashed] (0mm,48mm) -- (48mm,48mm);

\draw [dashed] (0mm,0mm) -- (0mm,48mm);
\draw [dashed] (8mm,0mm) -- (8mm,48mm);
\draw [dashed] (16mm,8mm) -- (16mm,48mm);
\draw [dashed] (24mm,16mm) -- (24mm,48mm);
\draw [dashed] (32mm,24mm) -- (32mm,48mm);
\draw [dashed] (40mm,32mm) -- (40mm,48mm);
\draw [dashed] (48mm,40mm) -- (48mm,48mm);

\node    at (4mm,44mm) {\footnotesize{$1-$$\frac{y_1}{y_2}$}};
\node    at (4mm,28mm) {\footnotesize{$1-$$\frac{y_3}{y_5}$}};
\node    at (4mm,20mm) {\footnotesize{$1-$$\frac{y_3}{y_7}$}};
\node    at (4mm,12mm) {\footnotesize{$1-$$\frac{y_3}{y_6}$}};

\node    at (12mm,28mm) {\footnotesize{$1-$$\frac{y_4}{y_5}$}};
\node    at (12mm,20mm) {\footnotesize{$1-$$\frac{y_4}{y_7}$}};
\node    at (12mm,12mm) {\footnotesize{$1-$$\frac{y_4}{y_6}$}};
\node    at (20mm,20mm) {\footnotesize{$1-$$\frac{y_6}{y_7}$}};

\draw[thick](0mm,40mm)--(8mm,40mm)--(8mm,48mm)--(0mm,48mm)--(0mm,40mm);

\draw [step=8mm,thick] (0mm,8mm) grid (16mm,32mm);

\draw[thick](16mm,16mm)--(24mm,16mm)--(24mm,24mm)--(16mm,24mm);

\end{tikzpicture}
\end{center}
\caption{Weights of squares of $D(w)$ for $w=2157634$.}\label{bot}
\end{figure}

To describe the Buch-Rim\'{a}nyi formula, consider the left-justified  array $\Delta_n$ with $n-i$ squares in row $i$.
Let $w=w(1)w(2)\cdots w(n)\in S_n$. For $1\leq i\leq n$, let
\[I(w,i)=\{w(j)\colon j>i,\ w(j)<w(i)\}\]
 be the set of entries in $w$ that are smaller than $w(i)$ but appear to the right of $w(i)$.
  Set $c(w,i)=|I(w,i)|$. It is clear that $0\leq c(w,i)\leq n-i$.
Let  $D(w)$ be the subset of $\Delta_n$
consisting of  the first $c(w,i)$ squares in the
  $i$-th row of $\Delta_n$, where $1\leq i\leq n$.
  Note that $D(w)$ corresponds to the
  bottom RC-graph of $w$, as defined by Bergeron and   Billey \cite{BB}.
  Assume that the values   in $I(w,i)$ are
\[w(j_1)<w(j_2)<\cdots<w(j_{c(w,i)}).\]
For a square $B\in D(w)$ in  row $i$
  and  column $k$,
equip $B$ with the  weight
\[\mathrm{wt}(B)=1-\frac{y_{w(j_k)}}{y_{w(i)}},\]
 see Figure \ref{bot} for an illustration.

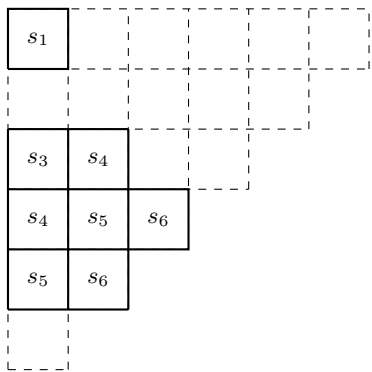
\begin{figure}[h]
\begin{center}
\begin{tikzpicture}(100,140)

\draw [dashed] (0mm,0mm) -- (8mm,0mm);
\draw [dashed] (0mm,8mm) -- (16mm,8mm);
\draw [dashed] (0mm,16mm) -- (24mm,16mm);
\draw [dashed] (0mm,24mm) -- (32mm,24mm);
\draw [dashed] (0mm,32mm) -- (40mm,32mm);
\draw [dashed] (0mm,40mm) -- (48mm,40mm);
\draw [dashed] (0mm,48mm) -- (48mm,48mm);

\draw [dashed] (0mm,0mm) -- (0mm,48mm);
\draw [dashed] (8mm,0mm) -- (8mm,48mm);
\draw [dashed] (16mm,8mm) -- (16mm,48mm);
\draw [dashed] (24mm,16mm) -- (24mm,48mm);
\draw [dashed] (32mm,24mm) -- (32mm,48mm);
\draw [dashed] (40mm,32mm) -- (40mm,48mm);
\draw [dashed] (48mm,40mm) -- (48mm,48mm);

\node    at (4mm,44mm) {\footnotesize{$s_1$}};
\node    at (4mm,28mm) {\footnotesize{$s_3$}};
\node    at (4mm,20mm) {\footnotesize{$s_4$}};
\node    at (4mm,12mm) {\footnotesize{$s_5$}};

\node    at (12mm,28mm) {\footnotesize{$s_4$}};
\node    at (12mm,20mm) {\footnotesize{$s_5$}};
\node    at (12mm,12mm) {\footnotesize{$s_6$}};
\node    at (20mm,20mm) {\footnotesize{$s_6$}};

\draw[thick](0mm,40mm)--(8mm,40mm)--(8mm,48mm)--(0mm,48mm)--(0mm,40mm);

\draw [step=8mm,thick] (0mm,8mm) grid (16mm,32mm);

\draw[thick](16mm,16mm)--(24mm,16mm)--(24mm,24mm)--(16mm,24mm);

\end{tikzpicture}
\end{center}
\caption{Labels of the squares of $D(w)$ for $w=2157634$.}\label{Lab}
\end{figure}

Given a subset  $D$ of $D(w)$, one can generate a word,
denoted $\mathrm{word}(D)$, as follows. Label the square of $D(w)$
in row   $i$  and column  $k$   by the simple transposition
  $s_{i+k-1}$, see Figure \ref{Lab} for an illustration.
  Then
$\mathrm{word}(D)$ is obtained by  reading off the labels of the squares
 in $D$ along the rows from top to bottom and right to left. For example,
 for the diagram $D=D(w)$ in Figure \ref{Lab}, we have
 \[\text{word}(D)=(s_1, s_4, s_3, s_6, s_5, s_4, s_6, s_5).\]

 A word $(s_{i_1}, s_{i_2}, \ldots, s_{i_m})$ is called a Hecke word of a permutation $u$ of length $m$ if
 \[(((s_{i_1}\ast s_{i_2})\ast s_{i_3})\ast\cdots)\ast s_{i_m}=u,\]
  where, for a permutation $w$, we define  $w\ast s_i$ to be $w$ if $\ell(w s_i)<\ell(w)$ and $w s_i$ otherwise.
  For example, $(s_1,s_2,s_1, s_2)$ is a Hecke word of $u=321$ of length 4 since
  \[((s_1\ast  s_2)\ast s_1)\ast s_2=((s_1s_2)\ast s_1)\ast s_2=(s_1s_2 s_1)\ast s_2=s_1s_2s_1=321.\]
  We note in passing that the
operation $\ast$ can be extended to an associative operation on the whole $S_n$;
this latter operation  is
the multiplication in the Hecke algebra associated to $S_n$ at $q=0$, see \cite[Chapter 7.4]{Hum}. Hence $\ast$   satisfies the associative property.  This means that the set of
permutations in $S_n$ forms  a monoid structure (0-Hecke monoid) under the operation $\ast$.

Write $\mathrm{Hecke}(D)=u$ if $\mathrm{word}(D)$ is
a Hecke word of a permutation  $u$. Notice that a Hecke word of $u$ of length $\ell(u)$ is a reduced word of $u$.
Note that for any $w\in S_n$, the
word $\mathrm{word}(D(w))$ is a reduced word of $w$, and therefore, if we multiply the letters of $\text{word}(D(w))$ using
either the $\ast$ product or the usual product of $S_n$, then we get $w$. That is, $\mathrm{Hecke}(D(w))=w$.

For any $u,v\in S_n$,  let
\[\mathcal{H}(u,v)=\{D\subseteq D(v)\,|\,\mathrm{Hecke}(D)=u\}.\]
For a   subset $D$ of $D(v)$, let
\begin{align}\label{wt}
\text{wt}(D)=\prod_{B \in D}\text{wt}(B).
\end{align}

\begin{theo}[\mdseries{Buch-Rim\'{a}nyi \cite[Theorem 2.1]{BR}}]\label{mt}
For permutations $u, v\in S_n$, we have
\begin{equation}\label{BR-m}
\G_u(y_v;y)=\sum_{D\in \mathcal{H}(u,v)}(-1)^{|D|-\ell(u)}\mathrm{wt}(D),
\end{equation}
where empty sums are interpreted as 0.
\end{theo}

We remark that in \cite{BR}, formula \eqref{BR-m} is
described in terms of the notation $C(\mathfrak{D}_v)$ and FK-graphs for $u$ with respect to $\mathfrak{D}_v$.
With the notation in this note, $D(v)$ can be obtained from  $C(\mathfrak{D}_v)$ by first reflecting along the main diagonal and then left-justifying the crossing positions. This operation also establishes a weight preserving bijection between the set $\mathcal{H}(u,v)$ and the set of FK-graphs for $u$ with respect to $\mathfrak{D}_v$.

\section{Elementary proof of Theorem \ref{mt}}

We need several  lemmas which follow directly  from the definition
of $\G_w(x;y)$.

\begin{lemma}\label{lm1}
Let $v=v's_i$ and $\ell(v)>\ell(v')$. If  $\ell(us_i)<\ell(u)$, then
\begin{equation}\label{fir}
\G_{u}(y_v;y)=\frac{y_{v'(i)}}{y_{v'(i+1)}}\G_{u}(y_{v'};y)+
\left(1-\frac{y_{v'(i)}}{y_{v'(i+1)}}\right)\G_{us_i}(y_{v'};y).
\end{equation}
\end{lemma}

\pf Applying \eqref{def} to $w=us_i$ and substituting $x_j$ with $y_{v'(j)}$, we have
\[
\G_{us_i}(y_{v'};y)
 =\frac{y_{v'(i)}\G_{u}(y_{v'};y)-y_{v'(i+1)}\G_{u}(y_v; y)}{y_{v'(i)}-y_{v'(i+1)}},
\]
which is equivalent to \eqref{fir}.
\qed

\begin{lemma}\label{lm2}
Let $v=v's_i$. If $\ell(us_i)>\ell(u)$, then
\begin{equation}\label{sec}
\G_{u}(y_v;y)=\G_{u}(y_{v'};y).
\end{equation}
\end{lemma}

\pf Applying \eqref{def} to $w=u$ and substituting $x_j$ with $y_{v(j)}$ and
$y_{v'(j)}$ respectively, we see that
\begin{align*}
\G_{u}(y_v;y)&=\frac{y_{v(i)}\G_{us_i}(y_v;y)-y_{v(i+1)}\G_{us_i}(y_{v'}; y)}{y_{v(i)}-y_{v(i+1)}},\\[5pt]
\G_{u}(y_{v'};y)&=\frac{y_{v'(i)}\G_{us_i}(y_{v'};y)-y_{v'(i+1)}\G_{us_i}(y_{v}; y)}{y_{v'(i)}-y_{v'(i+1)}},
\end{align*}
which, together with the fact that
 $v(i)=v'(i+1)$ and $v(i+1)=v'(i)$, implies  \eqref{sec}.
\qed

Let $\leq$ denote the (strong) Bruhat order on permutations of $S_n$. Recall that the Bruhat order is
the closure of the following covering relation: For $u, v\in S_n$, we say that $v$ covers
$u$  if there exists a transposition $t_{ij}$ such that $v=ut_{ij}$ and $\ell(v)=\ell(u)+1$. The following lemma is known, see \cite[Corollary 2.4]{BR} and the references therein.

\begin{lemma}\label{lm3}
We have $\G_u(y_v;y)=0$ whenever $u\not\leq v$ in the Bruhat order.
\end{lemma}

\pf  The idea in the proof of   \cite[(2.22)]{LLS}   for double Schubert polynomials
applies to   double Grothendieck polynomials, and we include a proof here  for the reader's  convenience.
Use descending induction on $\ell(u)$. The initial case is $u=w_0$.
Since $u\not\leq v$, we have $v\neq w_0$.
It is easily checked  from \eqref{DDD} that $\G_{w_0}(y_v; y)=0$.

We now consider the case $u\neq w_0$. Choose a position $i$ such that $u(i)<u(i+1)$.
Note that $u<us_i$.
Since $u\not\leq v$, we must have $us_i\not\leq v$.
We further claim that $us_i\not\leq vs_i$. This can be seen as follows. We have either $vs_i<v$ or $v<vs_i$ (depending on which of $\ell(vs_i)$ and  $\ell(v)$ is larger). If $vs_i<v$, then it is
clear that $us_i\not\leq vs_i$ since otherwise there would hold $u\leq v$.
It remains to verify the case $v<vs_i$. Suppose to the contrary that
$us_i\leq vs_i$. Then $u< vs_i$. Since $vs_i>v$ and $us_i>u$, applying
the Lifting Property  (see \cite[Proposition 2.2.7]{BBren}) to $u^{-1}$ and $(vs_i)^{-1}$, we obtain
that $u\leq v$, leading to a contradiction.
Now, by the definition in \eqref{def}  and by the induction hypothesis,
\[
\G_{u}(y_{v};y)
 =\frac{y_{v(i)}\G_{us_i}(y_{v};y)-y_{v(i+1)}\G_{us_i}(y_{vs_i}; y)}{y_{v(i)}-y_{v(i+1)}}=0,
\]
as desired.
\qed

\begin{lemma}\label{lm4}
Let $u\in S_n$ and $u'=us_i$ for some $i$ such that $\ell(us_i)<\ell(u)$.
Then,
\[
\G_u(y_u; y)=\left(1-\frac{y_{u(i+1)}}{y_{u(i)}}\right)\G_{u'}(y_{u'}; y).
\]
\end{lemma}

\pf
Apply Lemma \ref{lm1} to $v=u$ and $v'=u'$. The first addend on the right side
vanishes due to Lemma \ref{lm3}.
\qed

\begin{lemma}[\mdseries{Buch-Rim\'{a}nyi \cite[Corollary 2.6]{BR}}]\label{lm5}
For each $u\in S_n$, we have
\[
\G_u(y_u; y)=\prod_{i<j\atop u(i)>u(j)}\left(1-\frac{y_{u(j)}}{y_{u(i)}}\right).
\]
\end{lemma}

\pf Make descending induction on $\ell(u)$. The induction base for $u=w_0$
 is a  restatement of \eqref{DDD}. Assume that $u\neq w_0$.
 Then there exists some $1\le k<n$ such that $\ell(us_k)>\ell(u)$. Let $u'=us_k$. It is easy to see that the set
 \[\{(u'(i),u'(j))\,|\,i<j,\ u'(i)>u'(j)\}\]
is the union of the two disjoint sets
 \[\{(u(i),u(j))\,|\,i<j,\ u(i)>u(j)\}\cup \{(u(k),u(k+1))\}.\]
The proof follows by induction together with Lemma \ref{lm4}.
\qed

\vspace{10pt}

\noindent
{\it Proof of Theorem  \ref{mt}.}
The proof is by induction on $\ell(v)$.
 Let us first consider the case  $\ell(v)=0$, that is,
$v$ is the identity permutation $e$. If $u=e$, then it follows from Lemma \ref{lm5} (applied to $u=e$)
  that $\G_e(y_e;y)=1$. If $u\neq e$, then Lemma \ref{lm3} forces  that
$\G_u(y_e;y)=0$. So \eqref{BR-m} holds for $\ell(v)=0$.

Assume now  that $\ell(v)>0$. Let $s_r$ be the last
descent of $v$, that is,
$r$ is the largest index such that $v(r)>v(r+1)$.
Write $v=v's_r$.
Clearly, the bottom row of $D(v)$ lies in
 row $r$ of $\Delta_n$. The leftmost square in the bottom
row of $D(v)$, denoted $B_0$, has weight
\[\mathrm{wt}(B_0)=1-\frac{y_{v(r+1)}}{y_{v(r)}}=1-\frac{y_{v'(r)}}{y_{v'(r+1)}}.\]
Let  $u=u's_r$. There are two cases.

{\bf Case 1.}  $s_r$ is a descent of $u$.
By Lemma \ref{lm1} and by  induction hypothesis, we have
\begin{align}\label{FGX}
\G_u(y_v; y)
&=\frac{y_{v'(r)}}{y_{v'(r+1)}}\G_{u}(y_{v'};y)+
\left(1-\frac{y_{v'(r)}}{y_{v'(r+1)}}\right) \G_{u'}(y_{v'};y)\nonumber\\[5pt]
&=\left(1-\mathrm{wt}(B_0)\right) \sum_{D\in \mathcal{H}(u,v')} (-1)^{|D|-\ell(u)}\mathrm{wt}(D)+\mathrm{wt}(B_0)
\sum_{D\in \mathcal{H}(u',v')} (-1)^{|D|-\ell(u')} \mathrm{wt}(D) \nonumber\\[5pt]
&=\sum_{D\in \mathcal{H}(u,v')} (-1)^{|D|-\ell(u)}\mathrm{wt}(D)-\mathrm{wt}(B_0)\sum_{D\in \mathcal{H}(u,v')} (-1)^{|D|-\ell(u)}\mathrm{wt}(D)\nonumber\\[5pt]
&\quad\ \ +\mathrm{wt}(B_0)
\sum_{D\in \mathcal{H}(u',v')} (-1)^{|D|-\ell(u')} \mathrm{wt}(D).
\end{align}

To proceed, note that there is an obvious   bijection $\phi$ between  $D(v')$ and $D(v)\setminus \{B_0\}$.
 Since $s_r$ is the last descent of $v$,   we have $c(v',r)=0$,  $c(v',r+1)=c(v,r)-1$, and $c(v',i)=c(v,i)$ for $i\neq r, r+1$. Let $B\in D(v') $.
If $B$ lies  above row  $r$, then set $\phi(B)=B$.
Assume that $B$ lies in row  $r+1$ and column $j$,
then let $\phi(B)$ be the square of $D(v)\setminus \{B_0\}$
in row $r$ and column $j+1$.
By construction, $B$ and $\phi(B)$ are labeled by the same simple transposition. Moreover, it is easy to see that $\phi$ preserves the weight and   words, namely,
$\mathrm{wt}(B)=\mathrm{wt}(\phi(B))$ and $\text{word}(\phi(D))=\text{word}(D)$ for all $D\subseteq D(v')$.
Thus $\text{Hecke}(\phi(D))=\text{Hecke}(D)$ for all $D\subseteq D(v')$.

We claim that $\mathcal{H}(u,v)$ is the disjoint union of the following sets:
\begin{align*}
S_1&=\{\phi(D)\colon  D\in \mathcal{H}(u,v')\}, \\[5pt]
S_2&=\{\phi(D)\cup \{B_0\}\colon D\in \mathcal{H}(u,v')\}, \\[5pt]
S_3&=\{\phi(D)\cup \{B_0\}\colon D\in \mathcal{H}(u',v')\}.
\end{align*}
This can be easily seen as follows. Keep in mind that
$B_0$ is labeled by $s_r$.
Let $D\in \mathcal{H}(u,v)$. If $B_0\not\in D$, then $D\in S_1$. If $B_0\in D$, then
$\text{word}(D)$ is obtained from $\text{word}(D\backslash \{B_0\})$ by appending the letter $s_r$ at the end,
and thus we have $\text{Hecke}(D)=\text{Hecke}(D\backslash \{B_0\})\ast s_r$, and therefore  either $\mathrm{Hecke}(D\setminus \{B_0\})=u$ or $\mathrm{Hecke}(D\setminus \{B_0\})=u'$.
Hence either $D\in S_2$ or $D\in S_3$. Conversely, any $D\in S_1\cup S_2\cup S_3$ belongs to $\mathcal{H}(u,v)$, since $u\ast s_r=u'\ast s_r=u$.
By the above claim and in view of \eqref{FGX}, we obtain that
\begin{align*}
\G_u(y_v; y)
=\sum_{D\in S_1\cup S_2\cup S_3}(-1)^{|D|-\ell(u)}\mathrm{wt}(D)
=\sum_{D\in \mathcal{H}(u,v)}(-1)^{|D|-\ell(u)}\mathrm{wt}(D).
\end{align*}

{\bf Case 2.}   $s_r$ is not a descent of $u$.
Let $D\in \mathcal{H}(u,v)$.
We claim that $B_0\not\in D$. Suppose otherwise that $B_0 \in D$.
 Consider $D'=D\setminus \{B_0\}$.
If $s_r$ is a descent of $\mathrm{Hecke}(D')$, then $\mathrm{Hecke}(D)=\mathrm{Hecke}(D')$,
 while if $s_r$ is not a descent of $\mathrm{Hecke}(D')$, then $\mathrm{Hecke}(D)=\mathrm{Hecke}(D')\, s_r$.
 In both cases, $s_r$ is a descent of $u=\mathrm{Hecke}(D)$, leading to a
 contradiction.
Therefore, we see that $\mathcal{H}(u,v)=\{\phi(D)\,|\, D\in \mathcal{H}(u,v')\}$.
By Lemma \ref{lm2} and by induction hypothesis,
\begin{align*}
 \G_u(y_v; y) = \G_u(y_{v'}; y)
 =\sum_{D\in \mathcal{H}(u,v')} (-1)^{|D|-\ell(u)} \mathrm{wt}(D)
 =\sum_{D\in \mathcal{H}(u,v)} (-1)^{|D|-\ell(u)} \mathrm{wt}(D).
\end{align*}
This completes the proof.
\qed

\vskip 3mm \noindent {\bf Acknowledgments.}
We wish to thank the referees for valuable suggestions that greatly improve
the presentation of this note.
This work was supported by  the National Natural
Science Foundation of China (Grant No. 11971250).


\begin{thebibliography}{99}

\bibitem{BB}
N. Bergeron and S. Billey, RC-graphs and Schubert polynomials, Experiment. Math. 2 (4) (1993), 257--269.

\bibitem{BJS}
S. Billey, W. Jockusch and R.P. Stanley, Some combinatorial properties of Schubert
polynomials, J. Algebraic Combin. 2 (1993), 345--374.

\bibitem{BBren}
A. Bj\"orner and F. Brenti, Combinatorics of Coxeter Groups, Grad. Texts in Math.,
Vol. 231, Springer, New York, 2005.

\bibitem{BR}
A. Buch and R. Rim\'{a}nyi, Specializations of Grothendieck polynomials, C. R. Acad. Sci. Paris, Ser. I 339 (2004), 1--4.



\bibitem{FK1}
S. Fomin and A.N. Kirillov, Grothendieck polynomials and the Yang-Baxter equation, Proc. Formal Power Series and Alg. Comb. (1994), 183--190.

\bibitem{FK2}
S. Fomin and A.N. Kirillov, The Yang-Baxter equation, symmetric functions, and Schubert polynomials, in: Proceedings of the 5th Conference on Formal Power Series and Algebraic Combinatorics  (Florence, 1993),    Discrete Math. 153 (1996),   123--143.

\bibitem{FS}
S. Fomin and R.P. Stanley, Schubert polynomials and the NilCoxeter algebra, Adv.
Math. 103 (1994), 196--207.



\bibitem{Hum}
J.E. Humphreys,
Reflection Groups and Coxeter Groups,  Cambridge
Studies in Advanced Mathematics, No. 29, Cambridge Univ. Press,
Cambridge, 1990.

\bibitem{KM1}
A. Knutson and E. Miller, Gr\"obner geometry of Schubert polynomials, Ann. Math.
161 (2005), 1245--1318.

\bibitem{KM}
A. Knutson and E. Miller, Subword complexes in Coxeter groups, Adv. Math. 184 (2004), 161--176.

\bibitem{LLS}
T. Lam, S. Lee and M. Shimozono, Back stable Schubert calculus, arXiv:1806.11233v1.

\bibitem{LS}
 A. Lascoux and M.-P. Sch\"utzenberger,
Structure de Hopf de l'anneau de cohomologie et de l'anneau de Grothendieck d¡¯une vari\'et\'e
de drapeaux, C. R. Acad. Sci. Paris 295 (1982), 629--633.



\bibitem{LRS}
C. Lenart, S. Robinson  and F. Sottile, Grothendieck polynomials via permutation patterns and chains in the Bruhat order, Amer. J. Math. 128 (2006),   805--848.

\bibitem{Ma}
 I.G. Macdonald, Notes on Schubert Polynomials, Laboratoire de combinatoire et d'informatique math\'ematique (LACIM), Universit\'e du Qu\'ebec \'a Montr\'eal, Montreal, 1991.



\bibitem{YW}
A. Yong and A. Woo, A Gr\"{o}bner basis for the Kazhdan-Lusztig ideals, Amer. J. Math.  134 (2012), 1089--1137.

\end{thebibliography}
\end{document}